\documentclass[a4paper,12pt]{article}

\usepackage[latin1]{inputenc}
\usepackage[francais,english]{babel}
\usepackage{palatino}
\usepackage{fullpage}
\usepackage{url}
\usepackage{a4,amsmath,amssymb,amsfonts,amsthm}
\usepackage{amsrefs}

\newtheorem{theor}{Theorem}[section]
\newtheorem{prop}[theor]{Proposition}

\newtheorem{lemma}[theor]{Lemma}

\newtheorem{defin}[theor]{Definition}

\gdef\Bl{{\rm Bl}}
\gdef\Td{\op{\rm Td}}
\gdef\ch{\op{\rm ch}}

\gdef\CO{{\mathcal O}}

\def\M#1{\mathbb#1}     
\def\mR{\M{R}}           
\def\mZ{\M{Z}}           
\def\mN{\M{N}}           
\def\mQ{\M{Q}}       
\def\mC{\M{C}}

\gdef\beginProof{\par{\bf Proof. }}
\gdef\endProof{${\bf Q.E.D.}$\par}  
\gdef\ari#1{\op{\widehat{#1}}}
\gdef\mtr#1{\overline{#1}}
\gdef\c1{\op{{{\rm c}_1}}}
\gdef\ac1{\ari{\rm c}_1}
\gdef\aceq1{\ari{\rm c}_{{\rm eq},1}}

\gdef\Spec{\op{\rm Spec}}

\gdef\Vol{{\rm Vol}}

\gdef\refeq#1{(\ref{#1})}

\def\R{{\rm R}}

\def\Jac{{\rm Jac}}

\def\m2{{\mu_2}}
\def\rk{{\rm rk}}

\def\ul#1{\underline{#1}}

\def\op#1{\operatorname{#1}}

\def\be{\begin{equation}}
\def\ee{\end{equation}}
\def\BCOV{{\rm BCOV}}
\def\wt#1{\widetilde{#1}}
\def\top{{\rm top}}

\def\dR{{\rm dR}}
\def\Hdg{{\rm Hdg}}
\def\cl{{\rm cl}}

\author{D. R\"ossler \& V. Maillot}

\title{On a conjecture of H. Fang, Z. Lu and K.-I. Yoshikawa}

\begin{document}

\maketitle

\begin{abstract}
In \cite[Sec. 4, Conj. 4.17]{Fang-Lu-Yoshikawa}, Fang, Lu and Yoshikawa conjecture that 
a certain string-theoretic invariant of Calabi-Yau threefolds is a birational invariant. We prove 
a weak form of this conjecture.  
\end{abstract}

\bibliographystyle{plain}

\parindent=0pt


\section{Introduction}

Let $Y$ be a smooth projective variety of dimension $3$ over $\mC$. We suppose that $Y$ is a Calabi-Yau variety (in the restricted sense). By definition, this means that 
$H^1(Y,\CO_Y)=H^2(Y,\CO_Y)=0$ and that $\omega_Y:=\det(\Omega_Y)\simeq\CO_Y$. 

In \cite{Fang-Lu-Yoshikawa} (see also \cite[Sec. 2]{Yoshikawa-Analytic-torsion}), H. Fang, Z. Lu and K.-I. Yoshikawa  introduced the analytic invariant $\tau_\BCOV(Y(\mC))\in\mR^*_+$. See \cite[p. 177]{Fang-Lu-Yoshikawa} or 
Definition \ref{tbcovdef} below for the 
precise definition. They conjectured the following 
(see \cite[Sec. 4, Conj. 4.17]{Fang-Lu-Yoshikawa} and \cite[Sec. 2. Conj. 2.1]{Yoshikawa-Analytic-torsion}): 
if $Y$ and $Y'$ are birational Calabi-Yau varieties of dimension $3$ over $\mC$, then $\tau_\BCOV(Y(\mC))=\tau_\BCOV(Y'(\mC))$.\footnote{The conjecture made in \cite[Sec. 4, Conj. 4.17]{Fang-Lu-Yoshikawa} is only apparently 
weaker than the conjecture made in \cite[Sec. 2. Conj. 2.1]{Yoshikawa-Analytic-torsion}, because 
the topological types of $Y(\mC)$ and $Y'(\mC)$ coincide by a result of D. Huybrechts 
(see \cite[middle of p. 65]{Huybrechts-Compact}).}

H. Fang, Z. Lu and K.-I. Yoshikawa explain that their definition of $\tau_\BCOV$ is the mathematical formalisation 
of a definition made by the string-theorists M. Bershadsky, S. Cecotti, H. Ooguri and C.Vafa in 
\cite{Bershadsky-Holomorphic} and \cite{Bershadsky-Kodaira-Spencer} (the invariant $F_1(Y)$).
Their conjecture should be viewed as a "secondary" analog of the conjecture (which is now 
a theorem of Batyrev and Kontsevich; see \cite{Batyrev-Birational}) that the Hodge numbers of $Y(\mC)$ and $Y'(\mC)$ coincide. The latter conjecture was also motivated by physical considerations.

The purpose of this note is to describe the proof of the 
following arithmetic result, which is a step towards Yoshikawa's conjecture. 

Suppose that $X$ (resp. $X'$) is a smooth projective variety of dimension $3$ over $L$. Suppose that 
$X$ (resp. $X'$) is a Calabi-Yau variety (in the restricted sense). 

If $Z$ is a scheme, write as usual $D^b(Z):=D^b_c(Z)$ for the category derived from the homotopy classes of bounded 
complexes of coherent sheaves on $Z$. 

If $\sigma:L\hookrightarrow\mC$ is a subfield of $\mC$, we shall write 
$X_{\sigma}$ for the base change $X\times_{\Spec\ L,\sigma}\Spec\ \mC$ of 
$X$ to $\mC$ via the embedding $\sigma$. 

Now fix an embedding $\sigma:L\hookrightarrow\mC$. 

\begin{theor}
Let $T$ be a finite set of embeddings of $L$ into $\mC$. 

Let (S) be the statement : there exists $n\in\mN^*$ and $\alpha\in L^*$ such that for all $\tau\in T$, 
$$
\frac{\tau_\BCOV(X_\tau(\mC))}{\tau_\BCOV(X'_\tau(\mC))}=\sqrt[n]{|\tau(\alpha)|}.
$$

{\rm (A)} If $X_\sigma$ is birational to $X'_\sigma$ 
then (S) is verified.

{\rm (B)} If $D^b(X_\sigma)$ and $D^b(X'_\sigma)$ are equivalent as triangulated 
$\mC$-linear categories then (S) is verified.

\label{mainth}
\end{theor}

In particular, if $L=\mQ$ and  $X_\sigma$ is birational to $X'_\sigma$ 
 then there exists $n\in\mN^*$ such that $({\tau_\BCOV(X(\mC))}/{\tau_\BCOV(X'(\mC))})^n\in\mQ.$
 
Notice that by a theorem of Bridgeland (see \cite{Bridgeland-Flops}), (B) implies (A). 

We shall nevertheless give two separate proofs of (A) and (B).

Here is an outline of our proofs of (A) and (B). We first express the quantity 
$\tau_\BCOV$ in terms of arithmetic Chern numbers; this is made possible by the 
arithmetic Riemann-Roch theorem of Bismut-Gillet-Soulé \cite{Gillet-Soule-An-arithmetic}. 
To prove (A), we use the weak factorisation conjecture for birational maps (proved in \cite{Abramovich-Torification}) 
and some lemmas describing the effect of blow-up on some global Arakelov-theoretic invariants. 
To prove (B), we  make use of a theorem of Orlov, which asserts that if $D^b(X_\sigma)$ and $D^b(X'_\sigma)$ are equivalent as triangulated $\mC$-linear categories then $X_\sigma$ and $X'_\sigma$ are related 
by a Fourier-Mukai functor. We also use a theorem of Caldararu, which describes the effect of a Fourier-Mukai 
functor on the singular cohomology of $X_\sigma(\mC)$ and $X'_\sigma(\mC)$. 

{\bf Remark.} It is likely that Theorem \ref{mainth} is true without any 
restriction of finiteness on $T$. In particular, the quantity $\alpha$ should not depend on $T$. 
The reason for restricting the statement to finite $T$ is a (probably unnecessary) hypothesis of 
finiteness included in the definition of an arithmetic ring in Arakelov geometry. 









\section{The invariant $\tau_\BCOV$ and arithmetic Chern numbers on Calabi-Yau threefolds}
\label{arichernsec}

We shall apply the arithmetic Riemann-Roch theorem to certain 
vector bundles on $X$. In the following, we shall freely use the terminology of 
global Arakelov theory. For a concise summary of the necessary vocabulary, see \cite[Sec. 1]{Soule-Genres}. 

Let $f:X\to S:=\Spec\ L$ be the structure morphism. We may enlarge the set $T$ 
without changing the conclusion of Theorem \ref{mainth}, so we may assume that 
$T$ is conjugation-invariant. We view 
$L$ as an arithmetic ring, endowed with the set of embeddings 
$T$ into $\mC$. 
We endow $X(\mC):=\coprod_{\tau\in T}X_\tau(\mC)$ with  a conjugation-invariant Kähler form $\nu$. Let 
$\mtr{\Omega}:=\mtr{\Omega}_X$ be the sheaf of differentials of 
$X$, endowed with the metric induced by $\nu$. We write 
$\mtr{\omega}:=\mtr{\omega}_X$ for $\det(\mtr{\Omega}_X)$ and 
$\mtr{\Omega}^p:=\mtr{\Omega}^p_X$ for $\Lambda^p(\mtr{\Omega}_X)$. Furthermore, we shall 
write $H^q(Y,\mtr{\Omega}^p)$ for the $L$-vector space 
$R^q f_*(\Omega^p)$, endowed with the $L^2$-metric induced by $\nu$. 

Let $\mtr{\mathcal E}$ be the natural exact sequence of hermitian bundles
$$
0\to f^*f_*\mtr{\omega}\to\mtr{\omega}\to 0\to 0
$$
We let $\eta:=\wt{\ch}(\mtr{\mathcal E})$ be the Bott-Chern secondary class associated 
to $\mtr{\mathcal E}$, so that
$$
f^*f_*\mtr{\omega}-\mtr{\omega}=\eta
$$
in $\ari{K}_0(X)$, the arithmetic Grothendieck group of $X$. We shall write $\eta^0$ for 
the degree $0$ part of $\eta$. 

We apply the arithmetic Riemann-Roch theorem to $f$ and to the 
formal linear combination of hermitian bundles
$$
-\mtr{\Omega}^1+2\cdot\mtr{\Omega}^2-3\cdot\mtr{\Omega}^3.
$$
We obtain the equality 
\begin{eqnarray}
&&\ac1\big(\ \R^\bullet f_*[\ -\mtr{\Omega}^1+2\cdot\mtr{\Omega}^2-3\cdot\mtr{\Omega}^3\ ]\ \big)-a(\tau\big(\ -\mtr{\Omega}^1+2\cdot\mtr{\Omega}^2-3\cdot\mtr{\Omega}^3\ \big))\nonumber\\
&=&
f_*\big(\ari{\Td}(\mtr{\Omega}^\vee)\ari{\ch}\big(-\mtr{\Omega}^1+2\cdot\mtr{\Omega}^2-3\cdot\mtr{\Omega}^3\big)\big)^{(1)}\nonumber\\
&&-a(\large[\ \int_{X}R(\Omega)^\vee)\Td(\Omega^\vee)\ch(-{\Omega}^1+2\cdot{\Omega}^2-3\cdot{\Omega}^3)\ \large]^{(1)})\label{eqarr}
\end{eqnarray}
in $\ari{\rm CH}^1(L)_\mQ$, which is the first arithmetic Chow group of $L$, tensored with $\mQ$. 

Recall that $R(\cdot)$ is the $R$-genus of Gillet-Soulé and that $\tau(\cdot)$ is the Ray-Singer analytic 
torsion; see \cite[Introduction]{Gillet-Soule-An-arithmetic}.

We shall first analyse the various terms appearing in this equation. Write $\zeta_\mQ(s)$ for the 
evaluation of the Riemann zeta function at $s\in\mC$.

\begin{lemma} The equation
\begin{eqnarray*}
&&\ari{\Td}(\mtr{\Omega}^\vee)\ari{\ch}\big(-\mtr{\Omega}^1+2\cdot\mtr{\Omega}^2-3\cdot\mtr{\Omega}^3\big)\\
&=&-c^\top(\mtr{\Omega}^\vee)[\ a(\zeta_\mQ(0)\rk(\Omega))+\zeta_\mQ(-1)
\ac1(\mtr{\omega})+\textrm{terms of degree $>1$}\ ]\\
&=&
-c^\top(\mtr{\Omega}^\vee)[\ a(\zeta_\mQ(0)\rk(\Omega))+\zeta_\mQ(-1)
\ac1(f^*f_*\mtr{\omega})-a(\zeta_\mQ(-1)\eta^0)+\textrm{terms of degree $>1$}\ ]
\end{eqnarray*}
\label{lemgd}
holds in $\ari{\rm CH}^\bullet(X)_\mQ$. 
\end{lemma}
\beginProof
The proof is similar to the proof of \cite[Lemma 3.1]{Maillot-Rossler-On-the} so we omit it.
\endProof
Using the projection formula together with Lemma \ref{lemgd}, we obtain that
$$ 
\large[\ f_*(\ari{\Td}(\mtr{\Omega}^\vee)\ari{\ch}(-\mtr{\Omega}^1+2\cdot\mtr{\Omega}^2-3\cdot\mtr{\Omega}^3)\ \large]^{(1)}=
-a\big(\zeta_\mQ(-1)\ac1(f_*\mtr{\omega})\int_{X}c^\top(\Omega^\vee)\big)+
a\big(\zeta_\mQ(-1)\int_{X}c^\top(\mtr{\Omega}^\vee)\eta^0\big)
$$
in $\ari{\rm CH}^1(L)_\mQ$ (here the superscript $^{(1)}$ refers to part of degree $1$ in $\ari{\rm CH}^1(L)_\mQ$). 

We have the identity of cohomology classes 
$$
R(\Omega)^\vee)\Td(\Omega^\vee)\ch(-{\Omega}^1+2\cdot{\Omega}^2-3\cdot{\Omega}^3)=
-R(\Omega^\vee)c^\top(\Omega^\vee)[\ \zeta_\mQ(0)\rk(\Omega)+
\textrm{terms of degree $>0$}\ ]
$$
and so 
$$
\large[\ \int_{X}R(\Omega)^\vee)\Td(\Omega^\vee)\ch(\big(-{\Omega}^1+2\cdot{\Omega}^2-3\cdot{\Omega}^3\big)\ \large]^{(1)}=0
$$
since $R^1(\Omega^\vee)=0$ by assumption. 

As to the left-hand side of equation \refeq{eqarr}, we have 
\begin{eqnarray*}
&&\R^\bullet f_*[\ -\mtr{\Omega}^1+2\cdot\mtr{\Omega}^2-3\cdot\mtr{\Omega}^3\ ]\\
&=&
-[\ H^0(X,\mtr{\Omega}^1)-H^1(X,\mtr{\Omega}^1)+H^2(X,\mtr{\Omega}^1)\ ]+2[\ H^0(X,\mtr{\Omega}^2)-H^1(X,\mtr{\Omega}^2)+H^2(X,\mtr{\Omega}^2)\ ]\\
&-&3[\ H^0(X,\mtr{\Omega}^3)-H^1(X,\mtr{\Omega}^3)+H^2(X,\mtr{\Omega}^3)\ ]\\
&=&
-[\ -H^1(X,\mtr{\Omega}^1)+H^2(X,\mtr{\Omega}^1)\ ]
+2[\ -H^1(X,\mtr{\Omega}^2)+H^2(X,\mtr{\Omega}^2)\ ]
-3[\ H^0(X,\mtr{\Omega}^3)\ ]\\
&=&
-H^1(X,\mtr{\Omega}^1)-H^1(X,\mtr{\Omega}^2)-3\cdot H^0(X,\mtr{\CO})
\end{eqnarray*}
Here we used the fact that Serre duality is compatible with $L^2$-metrics (see \cite[p. 27, after eq. 9]{Gillet-Soule-Analytic}). 

Putting everything together, we get 
\begin{eqnarray*}
&&-\ac1(H^1(X,\mtr{\Omega}^1))-\ac1(H^1(X,\mtr{\Omega}^2))-3\ac1(H^0(X,\mtr{\CO}))-\tau\big(\ -\mtr{\Omega}^1+2\cdot\mtr{\Omega}^2-3\cdot\mtr{\Omega}^3\ \big)\\
&=&
-\zeta_\mQ(-1)\ac1(f_*\mtr{\omega})\int_{X}c^\top(\Omega^\vee)+
\zeta_\mQ(-1)\int_{X}c^\top(\mtr{\Omega}^\vee)\eta^0
\end{eqnarray*}
which implies that 
\begin{eqnarray*}
&&\ac1(-H^1(X,\mtr{\Omega}^1)-\ac1(H^1(X,\mtr{\Omega}^2)-\tau\big(\ -\mtr{\Omega}^1+2\cdot\mtr{\Omega}^2-3\cdot\mtr{\Omega}^3\ \big)\\
&=&
-\zeta_\mQ(-1)\ac1(f_*\mtr{\omega})\int_{X}c^\top(\Omega^\vee)+
\zeta_\mQ(-1)\int_{X}c^\top(\mtr{\Omega}^\vee)\eta^0-3\log\Vol(X(\mC),\nu)\\
&=&
{1\over 12}\ac1(f_*\mtr{\omega})\int_{X}c^\top(\Omega^\vee)-
{1\over 12}\int_{X}c^\top(\mtr{\Omega}^\vee)\eta^0-3\log\Vol(X(\mC),\nu)\\
\end{eqnarray*}
where 
$$
\Vol(X(\mC,\nu)):={1\over 3!(2\pi)^3}\int_{X(\mC)}\nu^3.
$$
Let us write
$$
\chi(X):=\int_{X}c^\top(\Omega^\vee).
$$
Note that $\chi(X)=\sum_{p,q}(-1)^{p+q}\dim_L(H^q(X,\Omega^p))$ by the generalized Gauss-Bonnet theorem. 
Notice also that $f_*\mtr{\omega}=H^0(X,\mtr{\Omega}^3)$ (by definition). Since 
$\zeta_\mQ(-1)=-1/12$, we see that 
\begin{eqnarray}
&&-\ac1(H^1(X,\mtr{\Omega}^2))-{1\over 12}\chi(X)\ac1(H^0(X,\mtr{\Omega}^3))\nonumber\\
&=&
-{1\over 12}\int_{X}c^\top(\mtr{\Omega}^\vee)\eta^0-3\log\Vol(X(\mC),\nu)
+\ac1(H^1(X,\mtr{\Omega}^1))\nonumber\\
&+&\tau\big(\ -\mtr{\Omega}^1+2\cdot\mtr{\Omega}^2-3\cdot\mtr{\Omega}^3\ \big)
\label{fundeq}
\end{eqnarray}
in $\ari{\rm CH}^1(L)$. 

The $L^2$-metric on $H^2(X(\mC),\mC)$ is induced from a 
Riemannian metric on the space $H^2(X(\mC),\mR)$. This is a consequence of 
the formula \cite[before Lemma 2.7]{Maillot-Rossler-On-the}. Let $\Vol_{L^2}(H^2(X(\mC),\mZ))$ be the volume of 
a fundamental domain of the lattice $H^2(X(\mC),\mZ)_{\rm free}$ in $H^2(X(\mC),\mR)$ for that 
metric. Here $H^2(X(\mC),\mZ)_{\rm free}$ is the largest direct summand of $H^2(X(\mC),\mZ)$, 
which is a free $\mZ$-module. 

\begin{lemma} The equality
$$
\ac1(H^1(X,\mtr{\Omega}^1))=a(-\bigoplus_{\tau\in T}\log(\Vol_{L^2}(H^2(X_\tau(\mC),\mZ))))
$$
holds in $\ari{\rm CH}^1(L)_\mQ$. 
\end{lemma}
\beginProof
Let $\tau\in T$ and let $e_1,\dots, e_r$ be a basis of $H^2(X_\tau(\mC),\mZ)_{\rm free}$. 
By definition, we have 
$$
\Vol_{L^2}(H^2(X_\tau(\mC),\mZ))=|e_1\wedge\dots\wedge e_r|^2
$$
where $|\cdot|$ refers to the natural norm on $\Lambda^r(H^2(X(\mC),\mC))$. Since 
$H^2(X_\tau(\mC),\mC)\simeq H^{1,1}(X(\mC))$ by hypothesis, we may conclude from 
the Lefschetz theorem on $(1,1)$-classes that the elements $e_i$ are classes 
of algebraic cycles $\bf e_i$ on $X_\tau$. 
Let $\ul{\tau}_0:K\hookrightarrow\mC$ be a field of definition for the $\bf e_i$, 
where $\ul{\tau}_0$ extends $\tau$. We may assume that $K$ is finite over $L$ (see 
\cite[proof of Prop. 1.5]{DMOS}). 
Write ${\bf e_i}^K$ for the model of ${\bf e_i}$ in $X_K$ and write $\cl_\dR$ for the cycle class 
map with values in de Rham cohomology. 
Let $\ul{\tau}:K\hookrightarrow \mC$ be another embedding of $K$ extending $\tau$.
Since by construction 
$$
H^2(X_{K,\ul{\tau}}(\mC),\mC)\simeq H^2_\dR(X_K/K)\otimes_{\ul{\tau}}\mC
$$
we see that  the elements 
$\cl_\dR({\bf e_i}^K)\otimes_{\ul{\tau}} 1$ form a basis of $H^2(X_{K,\ul{\tau}}(\mC),\mC)$. Furthermore, since 
$\cl_\dR({\bf e_i}^K)\otimes_{\ul{\tau}} 1=\cl_\dR({\bf e_i}^K\otimes_{\ul{\tau}} \mC)$, we see that 
the elements $\cl_\dR({\bf e_i}^K)\otimes_{\ul{\tau}} 1$ even form a basis of $H^2(X_{K,\ul{\tau}}(\mC),\mZ)_{\rm free}$. 
Furthermore, there is a natural identification 
$$
H^2(X_{K,\ul{\tau}}(\mC),\mZ)_{\rm free}\simeq H^2(X_{\tau}(\mC),\mZ)_{\rm free}
$$
which is an isometry for the $L^2$-metrics.

Now let $f:\Spec\ K\to\Spec\ L$ be the natural map. We view $\Spec\ K$ has an arithmetic variety over $\Spec\ L$. 
By the above, we have the equalities
\begin{eqnarray*}
\ac1(H^1(X_K,\mtr{\Omega}^1))&=&a(-2\bigoplus_{\ul{\tau}|\tau,\ \tau\in T}\log|\cl({\bf e_1}^K)\otimes_{\ul{\tau}} 1
\wedge\dots\wedge \cl({\bf e_r}^K)\otimes_{\ul{\tau}} 1|^2)\\
&=&a(-\bigoplus_{\ul{\tau}|\tau,\ \tau\in T}\log\Vol_{L^2}(X_\tau(\mC),\mZ))
\end{eqnarray*}
and thus 
\begin{eqnarray*}
[K:L]\ac1(H^1(X,\mtr{\Omega}^1))&=&f_*f^*\ac1(H^1(X,\mtr{\Omega}^1))=
f_*\ac1(H^1(X_K,\mtr{\Omega}^1))\\
&=&
[K:L]a(-\bigoplus_{\tau\in T}\log\Vol_{L^2}(X_\tau(\mC),\mZ))
\end{eqnarray*}
and we can conclude.
\endProof

The previous calculations motivate the following definition : 

\begin{defin}
\begin{eqnarray*}
\tau_\BCOV(X(\mC))&:=&\exp\large[\ -{1\over 12}\int_{X(\mC)}c^\top(\mtr{\Omega}^\vee_{X(\mC)})\eta^0-3\log\Vol(X(\mC),\nu)
-\log(\Vol_{L^2}(H^2(X(\mC),\mZ)))\\
&-&\tau(\mtr{\Omega}^1_{X(\mC)})+2\cdot\tau(\mtr{\Omega}^2_{X(\mC)})-3\cdot\tau(\mtr{\Omega}^3_{X(\mC)})\large\ ]
\end{eqnarray*}
\label{tbcovdef}
\end{defin}

It is proven in \cite[Sec. 4.4]{Fang-Lu-Yoshikawa} that $\tau_\BCOV(X)$ does not depend on the choice of $\nu$. 
Notice that 
equation \refeq{fundeq} together with the formula \cite[before Lemma 2.7]{Maillot-Rossler-On-the} already implies the weaker statement that 
$a(\tau_\BCOV(X))$ does not depend on $\nu$. 

The following equation summarizes the calculations made in this section:

\smallskip
\begin{equation}
\boxed{
\log(\tau_\BCOV(X(\mC)))=-\ac1(H^1(X,\mtr{\Omega}^2))-{1\over 12}\chi(X)\ac1(H^0(X,\mtr{\Omega}^3))
\textrm{\rm\ \ \ \ in\ \ $\ari{\rm CH}^1(L)_\mQ$}}
\label{fundeqfin}
\end{equation}

\section{Proof of Theorem \ref{mainth}}

With the equation \refeq{fundeqfin} in hand,  we see that Theorem \ref{mainth} is equivalent to 
the equation
\begin{equation}
-\ac1(H^1(X,\mtr{\Omega}^2))-{1\over 12}\chi(X)\ac1(H^0(X,\mtr{\Omega}^3))=
-\ac1(H^1(X',\mtr{\Omega}^2))-{1\over 12}\chi(X')\ac1(H^0(X',\mtr{\Omega}^3))
\label{tbp}
\end{equation}

\begin{lemma}
Let $L'$ be a finite field extension of $L$. We view $\Spec\ L'$ as an arithmetic variety over 
$L$. With this convention, the equation 
$$
-\ac1(H^1(X_{L'},\mtr{\Omega}^2))-{1\over 12}\chi(X_{L'})\ac1(H^0(X_{L'},\mtr{\Omega}^3))=
-\ac1(H^1(X'_{L'},\mtr{\Omega}^2))-{1\over 12}\chi(X'_{L'})\ac1(H^0(X'_{L'},\mtr{\Omega}^3))
$$
in $\ari{\rm CH}^1(L')_\mQ$ is equivalent to the equation \refeq{tbp}.
\label{intlem}
\end{lemma}
\beginProof
Let $f:\Spec\ L'\to\Spec\ L$ be the natural morphism. Using the projection formula, we compute
$$
[L':L]\ac1(H^1(X,\mtr{\Omega}^2))=f_*f^*\ac1(H^1(X,\mtr{\Omega}^2))=f_*\ac1(H^1(X_{L'},\mtr{\Omega}^2))
$$
and similarly
$$
[L':L]\ac1(H^0(X,\mtr{\Omega}^3))=f_*f^*\ac1(H^0(X,\mtr{\Omega}^3))=f_*\ac1(H^0(X_{L'},\mtr{\Omega}^3))
$$
If we combine these formulae with the analogous formulae for $X'$, we may conclude.
\endProof

Now notice that the group $\ari{\rm CH}^1(L')_\mQ$ (where $L'$ is viewed as an arithmetic variety over $L$) is naturally isomorphic to the homonymous group 
$\ari{\rm CH}^1(L')_\mQ:=\ari{\rm CH}^1(L',T')_\mQ$, which is the first arithmetic Grothendieck 
group of the arithmetic ring $L'$, endowed with the set
$$
T':=\{\tau':L'\hookrightarrow\mC|\tau'\in T\}
$$
of embeddings into $\mC$. Thus Lemma \ref{intlem} implies that the truth value of Theorem \ref{mainth} 
remains unchanged if we replace $L$ by a finite extension field $L'$ and $T$ by 
the set $T':=\{\tau':L'\hookrightarrow\mC|\tau'\in T\}$.

Before we begin with the proof, notice that by the formula \cite[before Lemma 2.7]{Maillot-Rossler-On-the}, the $L^2$-metric on 
$H^1(X,\mtr{\Omega}^2)$ is given by the formula
\begin{equation}
\langle \lambda,\kappa\rangle_{L^2}={i\over(2\pi)^3} \int_{X(\mC)}\lambda\wedge\mtr{\kappa}
\label{met12}
\end{equation}
and the $L^2$-metric on $H^0(X,\mtr{\Omega}^3)$ is given by the formula
\begin{equation}
\langle \lambda,\kappa\rangle_{L^2}={-i\over(2\pi)^3} \int_{X(\mC)}\lambda\wedge\mtr{\kappa}
\label{met03}
\end{equation}

In particular, these metrics do not depend on the choice of the K\"ahler form $\nu$. 

\subsection{Proof of (A)}
\label{fproofssec}

We now assume that there is a birational transformation from $X_\sigma$ to $X'_\sigma$.

\begin{lemma}
There is a birational transformation from $X_{\mtr{L}}$ to $X'_{\mtr{L}}$. 
\end{lemma}
\beginProof
This can be proven using a "spreading out" argument. We leave the details to the reader.
\endProof 

Notice that the birational transformation provided by the last Lemma has a model over a finite 
extension of $L$. Hence, by the discussion following Lemma \ref{intlem}, we may assume without loss
of generality that there is a birational transformation from $X$ to $X'$ defined over $L$. 

\begin{lemma}
$$
\ac1(H^0(X,\mtr{\Omega}^3))=\ac1(H^0(X',\mtr{\Omega}^3))
$$
\label{implem}
\end{lemma}
\beginProof
Let $\phi$ be a birational transformation from $X$ to $X'$. 
It is shown in \cite[Proof of Th. 8.19, chap. II]{Hartshorne-Algebraic} that there is an open set 
$U\subseteq X$ and a morphism $f:X\to X'$, with the following properties : $f$ induces $\phi$ and $\textrm{codimension}(U)\leqslant 2$. 
It is also shown in \cite[Proof of Th. 8.19, chap. II]{Hartshorne-Algebraic} that the maps
$$
H^0(X',\Omega^3)\stackrel{f^*}\longrightarrow H^0(U,\Omega^3)\stackrel{\textrm{restriction to $U$}}{\longleftarrow}
H^0(X,\Omega^3)
$$
are bijective. Thus, using the formula \refeq{met03}, we compute that
$$
\ac1(H^0(X',\mtr{\Omega}^3))=-\log|\int_{X'(\mC)}\lambda\wedge\mtr{\lambda}|=
-\log|\int_{X(\mC)}f^*(\lambda)\wedge\mtr{f^*(\lambda)}|=\ac1(H^0(X,\mtr{\Omega}^3))
$$
Here $\lambda\in H^0(X',\Omega^3)$ is any non-zero element.
\endProof

We recall the following theorem of Manin (and others). 

\begin{theor}
Let $Y$ be a smooth projective variety over $\mC$. Let 
$Z\hookrightarrow Y$ be a smooth closed subvariety of codimension $c$ of $Y$. 
Lety $\phi:\wt{Y}:=\Bl_{Z}(Y)\to Y$ be the blow-up of $Y$ along $Z$. Let $e:E\hookrightarrow \wt{Y}$ be the immersion of the exceptional 
divisor and let $\pi:E\to Z$ be the natural morphism. Let $\CO(1)$ be the tautological vector bundle 
on $E$. For any $k\in\mN$, there is an isomorphism of 
$\mQ$-Hodge structures
$$
H^k(Y,\mQ)\bigoplus\oplus_{l\geqslant 0}^{c-2}H^{k-2l}(Z(\mC),\mQ)(-l-1)\stackrel{\sim}{\to}H^k(\wt{Y},\mQ)
$$
given by the formula
$$
(\eta,\kappa_1,\dots,\kappa_{c-1})\mapsto (\phi^*\eta,e_*[\pi^*(\kappa_0)+\pi^*(\kappa_1)\cdot c_1(\CO(1))+
\pi^*(\kappa_2)\cdot c_1(\CO(1))^2+\dots+ \pi^*(\kappa_{c-2})\cdot c_1(\CO(1))^{c-2}])
$$
\label{thmanin}
\end{theor}
\beginProof See \cite{Manin-Correspondences}.\endProof 

\begin{lemma}
Let $C$ be a non-singular curve of genus $g$ over $L$. Then 
$$
\ac1(H^0(\Jac(C),\mtr{\Omega}^g))=\ac1(H^0(C,\mtr{\Omega}^1))+(g-1)\log(2\pi)
$$
in $\ari{\rm CH}^1(L)$, for any K\"ahler metrics on $C(\mC)$ and $\Jac(C)(\mC)$.
\label{faltlem}
\end{lemma}
\beginProof
See \cite[Exp. I, Lemme 3.2.1]{Szpiro-Mordell}. 
\endProof

\begin{prop}
Let $Y$ be a smooth projective threefold over $L$. Let 
$Z\hookrightarrow Y$ be a smooth closed subcurve of genus $g$ of $Y$. 
Let $\phi:\wt{Y}:=\Bl_{Z}(Y)\to Y$ be the blow-up of $Y$ along $Z$.  Then 
$$
\ac1(H^1(\wt{Y},\mtr{\Omega}^2))=\ac1(H^1({Y},\mtr{\Omega}^2))+\ac1(H^0(Z,\mtr{\Omega}^1))+2g\log(2\pi)
$$
for any K\"ahler metrics on $Y$, $Z$ and $\wt{Y}$. 
\label{propbl}
\end{prop}
\beginProof
Let $e:E\hookrightarrow \wt{Y}$ be the immersion of 
the exceptional divisor. Let $\pi:E\to Z$ be the natural morphism. 
By the Theorem  \ref{thmanin}, the map 
$$
H^1({Y},\Omega^2)\oplus H^0(Z,\Omega)\mapsto H^1(\wt{Y},\Omega^2)
$$
given by the formula
$$
(\eta,\kappa)\mapsto \phi^*(\eta)+e_*(\pi^*(\kappa))
$$
is an isomorphism. We compute 
\begin{eqnarray*}
&&{i\over (2\pi)^3}\int_{\wt{Y}(\mC)} \big(\phi^*(\eta_1)+e_*(\pi^*(\kappa_1))\big)\wedge 
\big(\phi^*(\mtr{\eta_2})+e_*(\pi^*(\mtr{\kappa_2}))\big)=\\
&=&
{i\over (2\pi)^3}\int_{\wt{Y}(\mC)}\phi^*(\eta_1)\wedge \phi^*(\mtr{\eta_2}) +
{i\over (2\pi)^3}\int_{\wt{Y}(\mC)}e_*\pi^*(z^*(\eta_1)\wedge\mtr{\kappa}_2)\\
&+&{i\over (2\pi)^3}\int_{\wt{Y}(\mC)}e_*\pi^*(z^*(\mtr{\eta}_2)\wedge\kappa_1)+
{i\over (2\pi)^3}\int_{\wt{Y}(\mC)}e_*(\pi^*(\kappa_1))\wedge e_*(\pi^*(\mtr{\kappa}_2))=\\
&=&
{i\over (2\pi)^3}\int_{\wt{Y}(\mC)}\phi^*(\eta_1)\wedge \phi^*(\mtr{\eta_2}) +
{i\over (2\pi)^3}\int_{\wt{Y}(\mC)}e_*(\pi^*(\kappa_1))\wedge e_*(\pi^*(\mtr{\kappa}_2))
\end{eqnarray*}
Now using the self-intersection formula (see for instance \cite[VI, 1., 1.4.2]{Fulton-Lang-Riemann-Roch}), we 
may compute 
\begin{eqnarray*}
&&{i\over (2\pi)^3}\int_{\wt{Y}(\mC)}e_*(\pi^*(\kappa_1))\wedge e_*(\pi^*(\mtr{\kappa}_2))=
{i\over (2\pi)^3}\int_{E(\mC)}e_*(\pi^*(\kappa_1)\wedge e^*e_*(\pi^*(\mtr{\kappa}_2)))=\\
&=&
{i\over (2\pi)^3}\int_{E(\mC)}e_*(c_1(\CO_E(-1))\wedge\pi^*(\kappa_1)\wedge \pi^*(\mtr{\kappa}_2))=\\
&=&
-{i\over (2\pi)^3}\int_{Z(\mC)}\kappa_1\wedge\mtr{\kappa}_2={1\over (2\pi)^2}\cdot {-i\over 2\pi}\int_{Z(\mC)}\kappa_1\wedge\mtr{\kappa}_2
\end{eqnarray*}
These formulae imply the conclusion of the proposition.
\endProof

\begin{lemma}
 Let $A$ and $B$ be abelian varieties over $\mtr{L}$ and let 
 $\phi:A_\sigma\to B_\sigma$ be an isogeny (over $\mC)$). Then there is 
 an isogeny $A\to B$ (over $\mtr{L}$).
\label{lemphi}
\end{lemma}
\beginProof
By spreading out. Left to the reader.
\endProof

\begin{lemma}
Let $A$ and $B$ be two abelian varieties over $L$ and suppose that 
there exists an isogeny $\phi:A\to B$ (over $L$). Suppose that $L$ contains a square 
root of $\deg(\phi)$. Then 
$$
\ac1(H^0(A,\mtr{\Omega}^g))=\ac1(H^0(B,\mtr{\Omega}^g))
$$
in $\ari{\rm CH}^1(L)$, for any choice of K\"ahler metrics on $A$ and $B$.
\label{isogflem}
\end{lemma}
\beginProof 
Let $\alpha_1,\dots,\alpha_g$ be a basis of the $L$-vector space 
$\Omega^g_{B}$. Using  the formula \cite[before Lemma 2.7]{Maillot-Rossler-On-the}, we see that 
for any embedding $\tau\in T$, we have
\begin{eqnarray*}
&&\langle \alpha_1\wedge\alpha_2\wedge\dots\wedge\alpha_g, \alpha_1\wedge\alpha_2\wedge\dots\wedge\alpha_g\rangle_{L^2}=\\
&=&\int_{B(\mC)}((\alpha_1\wedge\alpha_2\wedge\dots\wedge\alpha_g)\otimes_\tau 1)\wedge 
\mtr{((\alpha_1\wedge\alpha_2\wedge\dots\wedge\alpha_g)\otimes_\tau 1)}\\
&=&
\deg(\phi)^{-1}\int_{A(\mC)}((\phi^*(\alpha_1)\wedge\phi^*(\alpha_2)\wedge\dots\wedge\phi^*(\alpha_g))\otimes_\tau 1)\wedge 
\mtr{((\phi^*(\alpha_1)\wedge\phi^*(\alpha_2)\wedge\dots\wedge\phi^*(\alpha_g))\otimes_\tau 1)}
\end{eqnarray*}
and thus the mapping $H^0(B,\Omega^g)\to H^0(A,\Omega^g)$ given by 
the formula $\eta\mapsto(\sqrt{\deg(\phi)})\cdot\phi^*$ is an isometry of hermitian vector 
bundles. 
\endProof

Let now $\phi:X\dashrightarrow X'$ be a birational transformation. Let 
$X''$ be another smooth projective variety over $L$, together with morphisms 
$f:X''\to X$ and $g:X''\to X'$ such that $\phi\circ f$ and $g$ coincide as birational 
transformations. The variety $X''$ can be obtained as a desingularisation of 
the Zariski closure of the graph of $\phi$ in $X\times X'$. 

Denote by ${\mathcal PHS}(\mQ)$ the category of (pure) polarisable $\mQ$-Hodge structures.   

Using weak factorisation of birational maps (see \cite{Abramovich-Torification}) and Proposition \ref{propbl} and possibly 
replacing $L$ by one of its finite extensions, we conclude that there 
are curves $C_1,\dots, C_{r'}$ over 
${L}$ and numbers $s'_l\in\{-1,1\}$ so that
$$
H^3(X_\sigma(\mC),\mQ)+\sum_{l=1}^{r'}(-1)^{s'_l}H^1(C'_{l,\sigma},\mQ)(-1)=H^3(X''_\sigma(\mC),\mQ)
$$
in $K_0({\mathcal PHS}(\mQ))$ and so that 
$$
\ac1(H^1(X,\mtr{\Omega}^2))+\sum_{l=1}^{r'}(-1)^{s'_l}\ac1(H^0(C'_{l},\mtr{\Omega}^1))+2\sum_{l=1}^{r'}(-1)^{s'_l}
\textrm{genus}(C'_l)\log(2\pi)=\ac1(H^1(X'',\mtr{\Omega}^2)).
$$
Symmetrically,  there 
are curves $C''_1,\dots, C''_{r''}$ over 
${L}$ and numbers $s''_l\in\{-1,1\}$ so that
$$
H^3(X'_\sigma(\mC),\mQ)+\sum_{l=1}^{r''}(-1)^{s''_l}H^1(C''_{l,\sigma},\mQ)(-1)=H^3(X''_\sigma(\mC),\mQ)
$$
in $K_0({\mathcal PHS}(\mQ))$ and so that 
$$
\ac1(H^1(X',\mtr{\Omega}^2))+\sum_{l=1}^{r''}(-1)^{s''_l}\ac1(H^0(C''_l,\mtr{\Omega}^1))+2\sum_{l=1}^{r''}(-1)^{s''_l}
\textrm{genus}(C''_l)\log(2\pi)=\ac1(H^1(X'',\mtr{\Omega}^2)).
$$
Now by a theorem of Kontsevich (proved using motivic integration; see \cite{Looijenga-Motivic}) there is an isomorphism 
of $\mQ$-Hodge structures 
$
H^3(X_\sigma(\mC),\mQ)\simeq H^3(X'_\sigma(\mC),\mQ).
$ 
Thus 
$$
2\sum_{l=1}^{r'}(-1)^{s'_l}
\textrm{genus}(C'_l)\log(2\pi)= 2\sum_{l=1}^{r''}(-1)^{s''_l}
\textrm{genus}(C''_l)\log(2\pi).
$$
Furthermore, since the category of polarisable $\mQ$-Hodge structures is semi-simple,  there exists an isomorphism of 
$\mQ$-Hodge structures
$$
\bigoplus_{l,s'_l=1}H^1(\Jac(C'_l)_\sigma,\mQ)\bigoplus\bigoplus_{l,s''_l=-1}H^1(\Jac(C''_l)_\sigma,\mQ)\to\bigoplus_{l,s'_l=-1}H^1(\Jac(C'_l)_\sigma,\mQ)\bigoplus\bigoplus_{l,s''_l=1}H^1(\Jac(C''_l)_\sigma,\mQ)
$$
and thus an $\mtr{L}$-isogeny of abelian varieties
$$
\prod_{l,s'_l=1}\Jac(C'_l)\ \prod\ \prod_{l,s''_l=-1}\Jac(C''_l)\to\prod_{l,s'_l=-1}\Jac(C'_l)\ \prod\ \prod_{l,s''_l=1}\Jac(C''_l)
$$
Here we used Lemma \ref{lemphi}. Extend $L$ further so that the latter 
isogeny is defined over $L$. 
Then, by Lemma \ref{isogflem}, we have  
$$
\sum_{l}^{r'}(-1)^{s'_{l}}\ac1(H^0(\Jac(C'_{l}),\mtr{\Omega}^{\textrm{dim}\Jac(C'_{l})})=
\sum_{l}^{r''}(-1)^{s''_{l}}\ac1(H^0(\Jac(C''_{l}),\mtr{\Omega}^{\textrm{dim}\Jac(C''_{l})}).
$$
in $\ari{\rm CH}^1(L)_\mQ$. Using  Lemma \ref{faltlem}, we deduce that 
$$
\sum_{l}^{r'}(-1)^{s'_{l}}\ac1(H^0(C'_{l},\mtr{\Omega}^{1}))=
\sum_{l}^{r''}(-1)^{s''_{l}}\ac1(H^0(C''_{l},\mtr{\Omega}^{1})).
$$
so  that 
$$\ac1(H^1(X,\mtr{\Omega}^2))=\ac1(H^1(X',\mtr{\Omega}^2)).$$ 
Furthermore, by Lemma \ref{implem}, we have $$\ac1(H^0(X,\mtr{\Omega}^3))=\ac1(H^0(X',\mtr{\Omega}^3))$$
and by a theorem of Kontsevich (see \cite{Looijenga-Motivic}) we have $\chi(X)=\chi(X')$. This 
implies that 
$$
-\ac1(H^1(X,\mtr{\Omega}^2)-{1\over 12}\chi(X)\ac1(H^0(X,\mtr{\Omega}^3))=
-\ac1(H^1(X',\mtr{\Omega}^2)-{1\over 12}\chi(X')\ac1(H^0(X',\mtr{\Omega}^3)).
$$
Thus the equation \refeq{tbp} is verified and the theorem is proved.

\subsection{Proof of (B)}
\label{sproofssec}

We now assume that the categories 
$D^b(X_\sigma)$ and $D^b(X'_\sigma)$ are equivalent as triangulated 
$\mC$-linear categories. 

As a matter of notation, if $X_1\times X_2\times\dots\times X_t$ is a cartesian product of varieties, we shall write 
$$
\pi_{X_{i_1}X_{i_2}\dots X_{i_j}}^{X_1 X_2\dots X_t}:
X_1\times X_2\times\dots\times X_t\to X_{i_1}\times X_{i_2}\times\dots \times X_{i_j}
$$
for the natural projection. 

If $M$ (resp. $M'$) is an object in $D^b(X_\sigma)$ (resp. in $D^b(X'_\sigma)$, let $F_M$ (resp. $F_{M'}$) be the functor $D^b(X_\sigma)\to D^b(X'_\sigma)$ 
(resp. $D^b(X'_\sigma)\to D^b(X_\sigma)$) defined by the formula $$F_M(\cdot)=R^\bullet\pi^{X_\sigma X'_\sigma}_{X'_\sigma,*}(M\ul{\otimes} \pi^{X_\sigma X'_\sigma}_{X_\sigma,*}(\cdot))$$ 
(resp. $$F_{M'}(\cdot)=R^\bullet\pi^{X_\sigma X'_\sigma}_{X_\sigma,*}(M'\ul{\otimes} \pi^{X_\sigma X'_\sigma}_{X'_\sigma,*}(\cdot))$$). 
 The symbol $\ul{\otimes}$ refers to the derived 
tensor product and $R^\bullet f_*$ refers to the functor derived from the direct image functor. 

We shall make use of the following theorems. 

\begin{theor}[Orlov]
There exists an object $M$ (resp. $M'$) in $D^b(X_\sigma\times X'_\sigma)$ with the following properties. 
\begin{itemize}
\item[\rm (a)] The object 
$$
R^\bullet\pi_{X_\sigma X_\sigma ,*}^{X_\sigma X'_\sigma X_\sigma }(\pi^{X_\sigma X'_\sigma X_\sigma ,*}_{X_\sigma X'_\sigma }(M)\ul{\otimes}\pi^{X_\sigma X'_\sigma X_\sigma ,*}_{X'_\sigma X_\sigma }(M'))
$$
is isomorphic in $D^b(X_\sigma\times X_\sigma)$ to the image of the diagonal morphism in $X_\sigma\times X_\sigma$. 
\item[\rm (b)] The object 
$$
R^\bullet\pi_{X'_\sigma X'_\sigma ,*}^{X'_\sigma X_\sigma X'_\sigma }(\pi^{X'_\sigma X_\sigma X'_\sigma ,*}_{X'_\sigma X_\sigma }(M')\ul{\otimes}\pi^{X'_\sigma X_\sigma X'_\sigma ,*}_{X_\sigma X'_\sigma }(M))
$$
is isomorphic in $D^b(X'_\sigma\times X'_\sigma)$ to the image of the diagonal morphism in $X'_\sigma\times X'_\sigma$. 
\end{itemize}
\label{orlovth}
\end{theor}
\beginProof See \cite{Orlov-Derived}. \endProof

The last theorem is actually 
valid more generally if $X_\sigma$ (resp. $X'_\sigma$) is replaced by any smooth quasi-projective scheme over 
$\mC$ and if one assumes that $D^b(X_\sigma)$ and $D^b(X'_\sigma)$ are equivalent as triangulated 
$\mC$-linear categories. 

Write $\pi:X_\sigma\times X'_\sigma\to X_\sigma$ for the first projection and 
$\pi':X_\sigma\times X'_\sigma\to X'_\sigma$ for the second projection. 

\begin{theor}[Caldararu]
Let $M$ and $M'$ be objects satisfying the conditions 
(a) and (b) in Theorem \ref{orlovth}, then the map
$$
\Phi_M^H:H^\bullet(X_\sigma(\mC),\mQ)\to H^\bullet(X'_\sigma(\mC),\mQ)
$$
given by the formula
$$
\Phi_M^H(\beta):=\pi'_*(\pi^*(\beta)\cdot\sqrt{\Td}(X_\sigma\times X'_\sigma)\cdot\ch(M)))
$$
is an isomorphism and for any $k\in\mN$ we have 
$$
\Phi_M^H(\oplus_{p-q=k}H^{p,q}(X_\sigma(\mC)))=\oplus_{p-q=k}H^{p,q}(X'_\sigma(\mC))
$$
and furthermore, for any $\beta,\lambda\in H^3(X_\sigma(\mC),\mQ)$, we have
$$
\int_{X(\mC)}\beta\wedge\lambda=\int_{X'(\mC)}\Phi_M^H(\beta)\wedge\Phi_M^H(\lambda)
$$
\label{thcald}
\end{theor}
\beginProof
See \cite{Caldararu-The-Mukai} or \cite[5.2]{Huybrechts-Fourier-Mukai}.
\endProof

The last theorem is actually 
valid more generally if $X_\sigma$ (resp. $X'_\sigma$) is replaced by any smooth projective scheme of dimension $3$ over $\mC$.

Notice that Theorem \ref{thcald} implies that if its hypotheses are satisfied, then 
$$
\Phi_M^H(H^{2,1}(X_\sigma(\mC)))=H^{2,1}(X'_\sigma(\mC))
$$
and
$$
\Phi_M^H(H^{3,0}(X_\sigma(\mC)))=H^{3,0}(X'_\sigma(\mC)).
$$
Here we have used the fact that $X$ and $X'$ are Calabi-Yau varieties in 
the restricted sense. 

\begin{lemma}
There exists a finite field extension $K$ of $L$ and an object $M_0$ (resp. $M'_0$) of 
\mbox{$D^b(X_K\times X'_K)$} (resp. $D^b(X'_K\times X_K)$) such that 
\begin{itemize}
\item[$\textrm{\rm (a)}_K$] The object 
$$
R^\bullet\pi_{X_K X_K ,*}^{X_K X'_K X_K }(\pi^{X_K X'_K X_K ,*}_{X_K X'_K }(M)\ul{\otimes}\pi^{X_K X'_K X_K ,*}_{X'_K X_K }(M'))
$$
is isomorphic in $D^b(X_K\times X_K)$ to the image of the diagonal morphism in $X_K\times X_K$. 
\item[$\textrm{\rm (b)}_K$] The object 
$$
R^\bullet\pi_{X'_K X'_K ,*}^{X'_K X_K X'_K }(\pi^{X'_K X_K X'_K ,*}_{X'_K X_K }(M')\ul{\otimes}\pi^{X'_K X_K X'_K ,*}_{X_K X'_K }(M))
$$
is isomorphic in $D^b(X'_K\times X'_K)$ to the image of the diagonal morphism in $X'_K\times X'_K$. 
\end{itemize}
\label{desclem}
\end{lemma}
\beginProof 
 Let $\Delta:X_\sigma\hookrightarrow X_\sigma\times X_\sigma$ (
resp. $\Delta':X_\sigma\hookrightarrow X_\sigma\times X_\sigma$)  be the diagonal morphism. 
Let $U$ be a bounded complex of locally free sheaves on $X_\sigma\times X'_\sigma$ representing $M$ and let 
$U'$ be a bounded complex of locally sheaves on $X'_\sigma\times X_\sigma$ representing $M'$. 
Let $L_1$ be a finitely generated extension of $L$ (as a field), such that $U$ (resp. $U'$)
has a model over $X_{L_1}\times_{L_1}X_{L_1}$ (resp. $X'_{L_1}\times_{L_1}X'_{L_1}$). Let $S$ be an affine variety 
over $L$, which is smooth and irreducible and whose function field is isomorphic 
to $L_1$ as an $L$-algebra. After possibly replacing $S$ by one of its open affine subsets, we may find 
 bounded complexes of locally free sheaves $\wt{U}$ (resp. $\wt{U}'$) on $X_S\times_S X_S$ (resp. 
$X'_S\times_S X'_S$), which are models of $U$ and $U'$.

The conditions $\textrm{\rm (a)}$ and $\textrm{\rm (b)}$ in Theorem \ref{orlovth} are equivalent to the conditions : 
\begin{itemize}
\item There are isomorphisms of coherent sheaves
$$
R^0\pi_{X_\sigma X_\sigma ,*}^{X_\sigma X'_\sigma X_\sigma }(\pi^{X_\sigma X'_\sigma X_\sigma ,*}_{X_\sigma X'_\sigma }(M)\ul{\otimes}\pi^{X_\sigma X'_\sigma X_\sigma ,*}_{X'_\sigma X_\sigma }(M'))\simeq 
\Delta_*\CO_{X_\sigma}
$$
and 
$$
R^i\pi_{X_\sigma X_\sigma ,*}^{X_\sigma X'_\sigma X_\sigma }(\pi^{X_\sigma X'_\sigma X_\sigma ,*}_{X_\sigma X'_\sigma }(M)\ul{\otimes}\pi^{X_\sigma X'_\sigma X_\sigma ,*}_{X'_\sigma X_\sigma }(M'))\simeq 0
$$
for all $i\not= 0$;
\item  there are isomorphisms of coherent sheaves
$$
R^0\pi_{X'_\sigma X'_\sigma ,*}^{X'_\sigma X_\sigma X'_\sigma }(\pi^{X'_\sigma X_\sigma X'_\sigma ,*}_{X'_\sigma X_\sigma }(M')\ul{\otimes}\pi^{X'_\sigma X_\sigma X'_\sigma ,*}_{X_\sigma X'_\sigma }(M))\simeq 
\Delta_*\CO_{X'_\sigma}
$$
and
$$
R^i\pi_{X'_\sigma X'_\sigma ,*}^{X'_\sigma X_\sigma X'_\sigma }(\pi^{X'_\sigma X_\sigma X'_\sigma ,*}_{X'_\sigma X_\sigma }(M')\ul{\otimes}\pi^{X'_\sigma X_\sigma X'_\sigma ,*}_{X_\sigma X'_\sigma }(M))\simeq 0
$$
for all $i\not= 0$.
\end{itemize}
Thus, after possibly a further reduction of the size of $S$, we may assume that 
\begin{itemize}
\item there are isomorphisms of coherent sheaves
$$
R^0\pi_{X_S X_S ,*}^{X_S X_S' X_S }(\pi^{X_S X_S' X_S ,*}_{X_S X_S' }(\wt{U}){\otimes}\pi^{X_S X_S' X_S ,*}_{X_S' X_S }(\wt{U}'))\simeq 
\Delta_*\CO_{X_S}
$$
and 
$$
R^i\pi_{X_S X_S ,*}^{X_S X_S' X_S }(\pi^{X_S X_S' X_S ,*}_{X_S X_S' }(\wt{U}){\otimes}\pi^{X_S X_S' X_S ,*}_{X_S' X_S }(\wt{U}'))\simeq 0
$$
for all $i\not= 0$; 
\item there are isomorphisms of coherent sheaves
$$
R^0\pi_{X_S' X_S' ,*}^{X_S' X_S X_S' }(\pi^{X_S' X_S X_S' ,*}_{X_S' X_S }(\wt{U}'){\otimes}\pi^{X_S' X_S X_S' ,*}_{X_S X_S' }(\wt{U}))\simeq 
\Delta_*\CO_{X_S'}
$$
and
$$
R^i\pi_{X_S' X_S' ,*}^{X_S' X_S X_S' }(\pi^{X_S' X_S X_S' ,*}_{X_S' X_S }(\wt{U}'){\otimes}\pi^{X_S' X_S X_S' ,*}_{X_S X_S' }(\wt{U}))\simeq 0
$$
for all $i\not= 0$.
\end{itemize}
To see this, use the fact that the elements of the complexes $\wt{U}$ and $\wt{U}'$ are locally free and apply the theorem on cohomology and base-change (see \cite[chap. III, 7.7.4]{EGA}).

Now pick a closed point $s\in S$. The field $K:=\kappa(s)$ has all the properties we are looking for. 
\endProof

Now replace $L$ by a finite extension $K$ satisfying the conclusion of Lemma \ref{desclem}. 
Replace $T$ by the set $T_K$ of embeddings of $K$ into $\mC$ lying above embeddings 
in $T$. Recall that by Lemma \ref{intlem}, this does not restrict generality.

\begin{prop}
There are isometries of hermitian vector bundles
$$
H^1(X,\mtr{\Omega}^2)\simeq H^1(X',\mtr{\Omega}^2)
$$
and
$$
H^0(X,\mtr{\Omega}^3)\simeq H^0(X',\mtr{\Omega}^3)
$$
\label{isoprop}
\end{prop}
\beginProof
Let $M_0,M'_0$ be as provided by Lemma \ref{desclem}. 
Set $M_\sigma:=M_0\otimes_\sigma\mC$ and 
$M'_\sigma:=M'_0\otimes_\sigma\mC$. Since $\mC$ is flat as an $L$-algebra via $\sigma$, we see that $M_\sigma$ and $M'_\sigma$ satisfy properties (a) and (b) in Theorem \ref{orlovth}.

Furthermore, there are comparison isomorphisms 
$$
H^3(X_\sigma(\mC),\mC)\simeq (\bigoplus_{p+q=3}H^q(X,\Omega^p))\otimes_\sigma\mC
$$ and 
$$
H^3(X'_\sigma(\mC),\mC)\simeq (\bigoplus_{p+q=3}H^q(X,\Omega^p))\otimes_\sigma\mC,
$$
compatible with pull-backs, push-forwards and formation of Chern classes. 
We may thus conclude from Theorem \ref{thcald} that 
the morphism $\bigoplus_{p+q=3}H^q(X,\Omega^p)\to 
\bigoplus_{p+q=3}H^q(X,\Omega^p)$ 
given by the formula in Hodge cohomology 
$$
\Phi_M^{H_\Hdg}(\beta):=\pi'_*(\pi^*(\beta)\cdot\sqrt{\Td}(X\times X')\cdot\ch(M_0)))
$$
is an isomorphism. Therefore, again by Theorem \ref{thcald} (see remark after the theorem), the 
maps 
$$
\Phi_M^{H_\Hdg}\otimes_\tau\mC|_{H^1(X_\tau,\Omega^2)}:H^1(X_\tau,\Omega^2)\to H^1(X'_\tau,\Omega^2)$$
and
$$
\Phi_M^{H_\Hdg}\otimes_\tau\mC|_{H^0(X_\tau,\Omega^3)}:H^0(X_\tau,\Omega^3)\to H^0(X'_\tau,\Omega^3)
$$
are isometries for any $\tau\in T$. This implies the result. 
\endProof

We can now conclude the proof the Theorem \ref{mainth}. Indeed, by Proposition \ref{isoprop}, we 
have 
$$
\ac1(H^1(X,\mtr{\Omega}^2))=\ac1(H^1(X',\mtr{\Omega}^2))
$$
and
$$
\ac1(H^0(X,\mtr{\Omega}^3))=\ac1(H^0(X',\mtr{\Omega}^3))
$$
in $\ari{\rm CH}^1(L)_\mQ$. We conclude using equation \refeq{fundeq}.

\end{document}